\documentclass[12pt]{article}

\usepackage{amsmath,amssymb}
\def \a{{\mathfrak a}}
\def \adm{{\rm adm}}
\def \al{\alpha}
\def \aut{{\rm aut}}
\def \bs{\backslash}
\def \C{{\mathbb C}}
\def \CH{{\cal H}}
\def \CO{{\cal O}}

\def \CT{{\cal T}}

\def \End{{\rm End}}
\def \eps{{\varepsilon}}
\def \h{{\mathfrak h}}
\def \g{{\mathfrak g}}
\def \Ga{\Gamma}
\def \GL{{\rm GL}}
\def \Hom{{\rm Hom}}
\def \k{{\mathfrak k}}
\def \la{\lambda}
\def \L{{\mathfrak L}}

\def \Mat{{\rm Mat}}
\def \n{{\mathfrak n}}
\def \N{{\mathbb N}}
\def \p{{\mathfrak p}}
\def \PGL{{\rm PGL}}
\def \ph{{\varphi}}
\def \R{{\mathbb R}}

\def \SL{{\rm SL}}
\def \SO{{\rm SO}}
\def \st{{\rm st}}
\def \ua{{\underline a}}
\def \uan{{\underline {an}}}
\def \ul{\underline}
\def \un{{\underline n}}
\def \uk{{\underline k}}
\def \vol{{\rm vol}}
\def \Z{{\mathbb Z}}

\def \df{\ \begin{array}{c} _{\rm def}\\ ^{\displaystyle =}\end{array}\ }
\def \prf{{\bf Proof: }}
\def \qed{\ifmmode\eqno Q.E.D.
	\else\noproof\vskip 12pt plus 3pt minus 9pt \fi}
 	\def\noproof{{\unskip\nobreak\hfill\penalty50\hskip2em\hbox{}%
     \nobreak\hfill $\square$\parfillskip=0pt%
     \finalhyphendemerits=0\par}}

\def \({\left(}
\def \){\right)}
\def \={{\ =\ }}

\newcommand{\norm}[1]
	{\left|\hspace{-1pt}\left| #1\right|\hspace{-1pt}\right|}
\renewcommand{\sp}[2]{\left\langle #1,#2\right\rangle}

\newtheorem{theorem}{Theorem}[section]
\newtheorem{conjecture}[theorem]{Conjecture}
\newtheorem{lemma}[theorem]{Lemma}

\newtheorem{proposition}[theorem]{Proposition}

\begin{document}

\pagestyle{myheadings} \markright{INVARIANT TRIPLE PRODUCTS}

\title{Invariant triple products}
\author{Anton Deitmar}
\date{Int. J. Math. Math. Sci., Art. ID 48274, 22 pp. (2006).}

\maketitle

\tableofcontents

\section*{Introduction}

Let $G=\PGL_2(\R)$ and let $\pi_1,\pi_2,\pi_3$ be irreducible admissible smooth representations of $G$. Then the space of $G$-invariant trilinear forms on $\pi_1\times\pi_2\times\pi_3$ is at most one dimensional.
This has, in different contexts, been proved by
 Loke \cite{Lo}, Molchanov \cite{Mo} and Oksak \cite{Oksak}.
In this paper we ask for such a uniqueness result in the context of arbitrary semisimple groups.
We give evidence that for a given group $G$ uniqueness can only hold if $G$ is locally a product of hyperbolic groups.
For such groups we show uniqueness and for spherical vectors we compute the invariant triple products explicitly.

By a conjecture of Jacquet's, which has been proved in \cite{HK2}, triple products on $\GL_2$ are related to special values of automorphic $L$-functions, see also \cite{HK,HS,Jiang1,Jiang2}. The conjecture/theorem says that the existence of non-zero triple products is equivalent to the non-vanishing of the corresponding triple $L$-function at the centre of its functional equation.

The uniqueness of triple products in the $\PGL_2$-case mentioned above has been used in \cite{BR} to derive new bounds for automorphic $L^2$-coefficients.
This can also be done for higher dimensional hyperbolic groups, but, with the exception of the case treated in \cite{BR}, the results do not exceed those in \cite{KS}.
For completeness we include these computations in an appendix.

I thank J. Hilgert and A. Reznikov for helpful comments on the subject of this paper.

\section{Representations and Integral formulae}
Let $G$ be a connected semisimple Lie group with finite center.
Fix a maximal compact subgroup $K$.
Let $\hat G$ and $\hat K$ denote their unitary duals, i.e., the sets of isomorphism classes of irreducible unitary representations of $G$ resp. $K$.
Let $\pi$ be a continuous representation of $G$ on a locally convex topological vector space $V_\pi$.
Let $V_\pi'$ be the space of all continuous linear forms on $V_\pi$ and let $V_\pi^\infty$ be the space of \emph{smooth vectors}, i.e.,
$$
V_\pi^\infty\=\{ v\in V_\pi : x\mapsto \al(\pi(x)v) {\ \rm is\ smooth\ }\forall \al\in V_\pi'\}.
$$
The representation $\pi$ is called \emph{smooth} if $V_\pi=V_\pi^\infty$.

A representation $(\pi,V_\pi)$ is called \emph{admissible} if for each $\tau\in\hat K$ the space $\Hom_K(V_\tau,V_\pi)$ is finite dimensional.
Let $\hat G_{\adm}$ be the \emph{admissible dual}, i.e., the set of infinitesimal isomorphism classes of irreducible admissible representations.
A representation $\pi$ in $\hat G_\adm$ is called a \emph{class one} or \emph{spherical} representation if it contains $K$-invariant vectors.
In that case the space $V_\pi^K$ of $K$-invariant vectors is one-dimensional.
This is trivial for principal series representations (see below) and follows generally from Casselman's subrepresentation theorem which says that every $\pi\in\hat G_\adm$ is equivalent to a subrepresentation of a principal series representation.

The Iwasawa decomposition $G=ANK$ gives smooth maps
$$
\begin{array}{c}
\underline a \colon G\to A\\
\underline n \colon G\to N\\
\underline k \colon G\to K
\end{array}
$$
such that for every $x\in G$ one has $x=\ua(x)\un(x)\uk(x)$.
As an abbreviation we also define $\uan(x)=\ua(x)\un(x)$.
Let $\g_\R,\a_\R,\n_\R,\k_\R$ denote the Lie algebras of $G,A,N,K$ and let $\g,\a,\n,\k$ be their complexifications.

For $x\in G$ and $k\in K$ we define
$$
k^x\df \uk(kx).
$$

\begin{lemma}
The rule $k\mapsto k^x$ defines a smooth (right) group action of $G$ on $K$.
\end{lemma}

\prf
The map $\underline k$ gives a diffeomorphism $AN\bs G\to K$ and the action under consideration is just the natural right action of $G$ on $AN\bs G$.
\qed

\begin{lemma}
For $f\in C(K)$ and $y\in G$ we have the integral formula
$$
\int_{K} f(k)\, dk\= \int_{K} \ua(ky)^{2\rho}\, f(k^y)\, dk,
$$
or,
$$
\int_{K} f(k^y)\, dk\= \int_{K} \ua(ky^{-1})^{2\rho}\, f(k)\, dk.
$$
Here $\rho\in\a^*$ is the modular shift, i.e., $a^{2\rho}=\det(a|\n)$.
\end{lemma}

\prf
The Iwasawa integral formula implies
$$
\int_G g(x)\, dx\= \int_{AN}\int_K g(ank)\, dk\, dan.
$$
Let $f\in C(K)$ and choose  a function $\eta\in C_c(AN)$ such that $\eta\ge 0$ and $\int_{AN}\eta(an)\, dadn=1$.
Let $g(x)=\eta(\uan(x))f(\uk(x)$.
Then $g\in C_c(G)$ and
\begin{eqnarray*}
\int_G g(x)\, dx &=& \int_{AN}\eta(an)\, dan\, \int_K f(k)\, dk\\
&=& \int_Kf(k)\, dk.
\end{eqnarray*}
On the other hand, since $G$ is unimodular, this also equals
\begin{eqnarray*}
\int_G g(xy)\, dx &=& \int_G\eta(\uan(xy))f(\uk(xy))\, dx\\
&=& \int_{AN}\int_K \eta(\uan(anky))f(\uk(ky))\, dank\\
&=& \int_K\( \int_{AN} \eta(\uan(anky))\, dan\) f(k^y)\, dk\\
&=& \int_K\( \int_{AN} \eta(an\uan(ky))\, dan\) f(k^y)\, dk\\
&=& \int_K \ua(ky)^{2\rho}\( \int_{AN} \eta(an)\, dan\) f(k^y)\, dk\\
&=& \int_K \ua(ky)^{2\rho} f(k^y)\, dk.
\end{eqnarray*}

The second assertion follows from the first by replacing $f$ with
$\tilde f(k)\= f(k^y)$ and then $y$ with $y^{-1}$.
\qed

Let $M$ be the centraliser of $A$ intersected with $K$, then $P=MAN$ is a minimal parabolic subgroup of $G$.
The inclusion map $K\hookrightarrow G$ induces a diffeomorphism $M\bs K \to P\bs G$ and in this way we get a smooth $G$-action on $M\bs K$.
An inspection shows that this action is given by
$Mk\mapsto Mk^x$ for $k\in K$, $x\in G$.

\section{Trilinear products}

Let $\pi_1,\pi_2,\pi_3$ be three admissible smooth representations of the group $G$ and let $\CT\colon V_{\pi_1}\times V_{\pi_2}\times V_{\pi_3}\to \C$ be a continuous $G$-invariant trilinear form, i.e.,
$$
\CT(\pi_1(x)v_1, \pi_2(x)v_2, \pi_3(x)v_3)\= \CT(v_1,v_2,v_3)
$$
for all $v_j\in V_{\pi_j}$ and every $x\in G$.

We want to understand the space of all trilinear forms $\CT$ as above.
In this paper we will only consider principal series representations, the general case will be considered later.
So we assume that $\pi_1,\pi_2,\pi_3$ are principal series representations. This means that there are given a minimal parabolic $P=MAN$, irredicible representations $\sigma_j\in\hat M$, and $\la_j\in\a^*$ for $j=1,2,3$.
Each pair $(\sigma_j,\la_j)$ induces a continuous group homomorphism $P\to\GL(V_{\sigma_j})$ by $man\mapsto a^{\la_j+\rho}\sigma(m)$, which in turn defines a $G$-homogeneous vector bundle $E_{\sigma_j,\la_j}$ over $P\bs G$.
The representation $\pi_j$ is the $G$-representation on the space of smooth sections $\Ga^\infty(E_{\sigma_j,\la_j})$ of that bundle.
In other words, $\pi_j$ lives on the space of all $C^\infty$ functions $f\colon G\to V_{\sigma_j}$ with
$$
f(manx)\= a^{\la_j+\rho} \sigma_j(m)f(x)
$$
for all $m\in M$, $a\in A$, $n\in N$, $x\in G$.
The representation $\pi_j$ is defined by $\pi_j(y)f(x)=f(xy)$.
Every such $f$ is uniquely determined by its restriction to $K$ which satisfies $f(mk)=\sigma_j(m)f(k)$, i.e., $f$ is a section of the $K$-homogeneous bundle $E_{\sigma_j}$ on $M\bs K$ induced by $\sigma_j$. So the representation space can be identified with
$V_{\pi_j}\cong \Ga^\infty(E_{\sigma_j})$.
Thus a trilinear form $\CT$ is a distribution on the vector bundle $E_\sigma=E_{\sigma_1}\boxtimes E_{\sigma_2}\boxtimes E_{\sigma_3}$ over $M\bs K\times M\bs K\times M\bs K$.
Here $\boxtimes$ denotes the outer tensor product.

For $f_1,f_2,f_3\in C^\infty(M\bs K)$ we write
$\CT(f_1,f_2,f_3)$ for the expression
$$
\int_{M\bs K\times M\bs K\times M\bs K} \phi(k_1,k_2,k_3)\,\left[f_1(k_1)\boxtimes f_2(k_2)\boxtimes f_3(k_3)\right]\, dk_1\, dk_2\, dk_3,
$$
where $\phi$ is the kernel of $\CT$.

The group $G$ is called a \emph{real hyperbolic group} if it is locally isomorphic to ${\rm SO}(d,1)$ for some $d\ge 2$.

On $Y=(P\bs G)^3$ we consider the $G^3$-homogeneous vector bundle $E_{\sigma,\la}$ given by $E_{\sigma,\la}=E_{\sigma_1,\la_1}\boxtimes E_{\sigma_2,\la_2}\boxtimes E_{\sigma_3,\la_3}$.
Next $Y$ can be viewed as a $G$-space via the diagonal action and so $E_{\sigma,\la}$ becomes a $G$-homogeneous line bundle on $Y$.

We are going to impose the following condition on the induction parameters $\la_1,\la_2,\la_3$.
We assume that $\sum_{j=1}^3 \eps_j(\la_j+\rho)\ne 0$ for any choice of $\eps_j\in\{ \pm 1\}$.
In other words this means
\begin{itemize}
\item $\la_1+\la_2+\la_3+ 3\rho\ne 0$,
\item $\la_1+\la_2-\la_3 +\rho\ne 0$,
\item$\la_1-\la_2 +\la_3 +\rho\ne 0$, and
\item$-\la_1+\la_2 +\la_3 +\rho\ne 0$.
\end{itemize}

\begin{theorem}
Assume the parameters $\la_1,\la_2, \la_3$ satisfy the above condition.
Let $Y$ be the $G$-space $(P\bs G)^3$.
If there is an open $G$-orbit in $Y$, then the dimension of the space of invariant trilinear forms on smooth principal series representations is less than or equal to
$$
\sum_o \dim(\sigma_1\otimes\sigma_2\otimes\sigma_3)^{M_o},
$$
where the sum runs over all open orbits $o$ and $M_o$ is the stabilizer group of a point in the orbit $o$ which is chosen so that $M_o$ is a subgroup of $M$.
If all induction parameters are imaginary (unitary induction), then there is equality.

In particular, if $\pi_1,\pi_2,\pi_3$ are class one representations, then the dimension is less than or equal to the number of open orbits in $Y$.

There is an open orbit if and only if $G$ is locally isomorphic to a product of hyperbolic groups.

For $G=\SO(2,1)^0$ the number of open orbits is $2$, for $G=\SO(2,1)$ or $G=\SO(d,1)^0$, $d\ge 3$, the number is $1$.
Here $\SO(d,1)^0$ is the connected component of the Lie group $\SO(d,1)$.
\end{theorem}

The Proof is based on the following lemma.

\begin{lemma}
Let $G$ be a Lie group and $H$ a closed subgroup.
Let $X=G/H$ and let $E\to X$ be a smooth $G$-homogeneous vector bundle.
Let $\CT$ be a distribution on $E$, i.e. a continuous linear form on $\Ga_c^\infty(E)$.
Suppose that $\CT$ is $G$-invariant, i.e., $\CT(g.s)=\CT(s)$ for every $s\in \Ga_c^\infty(E)$.
Then $\CT$ is given by a smooth $G$-invariant section of the dual bundle $E^*$.

Let $(\sigma, V_\sigma)$ be the representation of $H$ on the fibre $E_{eH}$ and let $(\sigma^*,V_{\sigma^*})$ be its dual.
Then the space of all $G$-invariant distributions on $E$ has dimension equal to the dimension $\dim V_{\sigma^*}^H$ of $H$-invariants.
So in particular, if $\sigma$ is irreducible, this dimension is zero unless $\sigma$ is trivial, in which case the dimension is one.
\end{lemma}

\prf
The equation $\CT(g.s)=\CT(s)$, i.e., $g.\CT=\CT\ \forall g\in G$ implies $X.\CT=0$ for every $X\in \g_\R$, the real Lie algebra of $G$.
Let $\h_\R$ be the Lie algebra of $H$ and choose a complementary space $\p_\R$ for $\h_\R$ such that $\g_\R=\h_\R\oplus \p_\R$.
Let $X_1,\dots,X_n$ be a basis of $V$ and let
$$
D\= X_1^2+X_2^2+\cdots +X_n^2\ \in\ {\rm U}(\g_\R).
$$
We show that $D$ induces an elliptic differential operator on $E$.
By $G$-homogeneity, it suffices to show this at a single point.
So let  $P=\exp(\p_\R)$.
In a neighbourhood $U$ of the unit in $G$ there are smooth maps $\ul h:U\to H$ and $\ul p :U\to P$ such that $\ul h(x)\ul p(x)=x$ for $x\in U$.

The sections of $E$ can be identified with the smooth maps $s:G\to V_\sigma$ with $s(hx)=\sigma(h)s(x)$ for $h\in H$ and $x\in G$.
We can attach to each section $s$ a map $f_s$ on $\p$ with values in $V_\sigma$ by $f_s(Y)=s(\exp(Y))$.
The action of $X\in\p_\R$ on the section $s$ is described by
\begin{eqnarray*}
f_{Xs}(Y) &=& \left.\frac d{dt}\right|_{t=0} s(\exp(Y)\exp(tX))\\
&=& \left.\frac d{dt}\right|_{t=0} \sigma(\ul h(\exp(Y)\exp(tX))\, s(\ul p(\exp(Y)\exp((tX)).
\end{eqnarray*}
Let $A_X(Y)=\left.\frac d{dt}\right|_{t=0}\sigma(\ul h(\exp(Y)\exp(tX))\in\End(V_\sigma)$.
Then the Leibniz rule implies
$$
f_{Xs}(Y)\= A_X(Y) f_s(Y) + Xf_s(Y).
$$
The first summand is of order zero and the second of order one.
Moreover, the second summand at $Y=0$ coincides with the coordinate-derivative in the direction of $X$.
This implies that the leading symbol of $D$ at $eH$ is $\xi_1^2+\dots +\xi_n^2$ and so $D$ is elliptic.
The distributional equation $D\CT=0$ then implies that $\CT$ is given by a smooth section.

For the second assertion of the lemma recall that a $G$-invariant section is uniquely determined by its restiction to the point $eH$ which must be invariant under $H$.
\qed

For the proof of the theorem we will need to investigate the $G$-orbit structure of $Y=(P\bs G)^3$.
First note that since the map $Mk\mapsto Pk$ is a $K$-isomorphism from $M\bs K$ to $P\bs G$, the $K$-orbit of every $y\in Y$ contains an element of the form $(y_1,y_2,1)$.
Hence the $P=MAN$-orbit structure of $(P\bs G)^2$ is the same as the $G$-orbit structure of $Y$.
By the Bruhat decomposition, the $P$-orbits in $P\bs G$ are parametrized by the Weyl group $W=W(G,A)$, where the unique open orbit is given by $Pw_0P$, here $w_0$ is the long element of the Weyl group.
Note that the $P$-stabilizer of $Pw_0\in P\bs G$ equals $AM$.
This implies that the $G$-orbits in $Y$ of maximal dimension are in bijection to the $AM$-orbits in $P\bs G$ of maximal dimension via the map $PxAM\mapsto (x,w_0,1).G$.
Again by Bruhat decomposition it follows that the latter are contained in the open cell $Pw_0P=Pw_0N$.
So the $G$-orbits of maximal dimension in $Y$ are in bijection to the $AM$-orbits in $N$ of maximal dimension, where $AM$ acts via the adjoint action.
The exponential map $\exp\colon \n_\R\to N$ is an $AM$-equivariant bijection, so we are finally looking for the $AM$-orbit structure of the linear adjoint action on $\n_\R$.

We will now prove that
there is an open orbit if and only if $G$ is locally a product of real hyperbolic groups.
So suppose that $Y$ contains an open orbit.
Then $\n_\R$ contains an open $AM$-orbit, say $AM.X_0$.
Let $\phi^+$ be the set of all positive restricted roots on $\a={\rm Lie}(A)$.
Decompose $\n_\R$ into the root spaces
$$
\n_\R\=\bigoplus_{\al\in\phi^+} \n_{\R,\al}.
$$
On each $\n_{\R,\al}$ install an $M$-invariant norm $\norm{.}_\al$.
This is possible since $M$ is compact.
Consider the map
\begin{eqnarray*}
\psi\colon \n_\R & \to & \prod_{\al\in\phi^+} \R\\
x& \mapsto & \prod_\al \norm{x}_\al.
\end{eqnarray*}
Since the orbit $AM.X_0$ is open, the image $\psi(AM.X_0)$ of the orbit must contain a nonempty open set.
Away from the set $\{ X\in\n_\R : \exists \al: \norm{x}_\al=0\}$ the map $\psi$ can be chosen differentiable.
Since the norms are invariant under $M$, one gets a smooth map
\begin{eqnarray*}
A &\to & \R^{|\phi^+|}\\
a & \mapsto & \psi(a.X_0),
\end{eqnarray*}
whose image contains an open set.
This can only happen if the dimension of $A$ is at least as big as $|\phi^+|$ and the latter implies that $G$ is locally a product of real rank one groups.
Now by Araki's table (\cite{Helg}, pp. 532-534) one knows that these real rank one groups must all be hyperbolic, because otherwise there would be two different root lengths.

For the converse direction let $G$ be locally isomorphic to ${\rm SO}(d,1)$.
We have to show that there is an open $AM$-orbit in $\n_\R$.
This, however, is clear as the action of $AM$ on $\n_\R$ is the natural action of $\R^\times_+\, \times\, {\rm SO}(d-1)$ on $\R^{d-1}$, hence there are two orbits, the zero orbit and one open orbit.

We will now show that if there is an open orbit, then there are no invariant distributions supported on lower dimensional orbits.
For this it suffices to consider the case $G={\rm SO}^0(d,1)$.
For simplicity, we only consider the trivial bundle, i.e., functions instead of sections. The general case is similar.
The orbit structure is as follows.
For $d=2$ the group $M=\SO(d-1)$ is trivial and so there are two open $AM$-orbits in $N\cong\R$ in this case. If $d>2$ there is only one open orbit.
Since the proof is very similar in the case $d=2$ we will now retrict ourselves to the case $d>2$.
The open orbit $[w_0n_0,w_0,1]$, given by some $n_0\in N$, contains in its closure the orbits $[w_0,w_0,1]$ and $[w_0,1,1]$ which in turn contain in their closures the orbit $[1,1,1]$.
$$
\begin{array}{rcl}
&[w_0n_0,w_0,1]& \\
&&\\
\cup &\cup & \cup \\
&&\\
\ [w_0,w_0,1] & [w_0,1,w_0] & [w_0,1,1]\\
&&\\
\cup &\cup & \cup\\
&&\\
&[1,1,1]&
\end{array}
$$
Let $\CO=[w_0n_0,w_0,1]$ be the open orbit.
Let $x_1=(w_0,w_0,1)$, then one has $[w_0,w_0,1]=x_1.G$.
Consider $A\cong\R_+^\times$ as a subset of $\R$ suitably normalised, then one can write
$$
x_1\=\lim_{a\to 0} x_0.am.
$$
Let $\CT$ be a $G$-invariant distribution supported on the closure of the orbit of $x_1$.
Since $\CT$ is $G$-invariant, it satisfies $X.\CT=0$ for everey $X\in\g$. 
Hence the wave front set $WF(\CT)\subset T^*Y$ is a $G$-invariant subset of the normal bundle of the manifold $x_1.G$.
This implies that $\CT$ is of order zero along the manifold $x_1.G$.
By the $G$-invariance it follows that $\CT$ is of the form
$$
\CT(f)\=\int_{x_1.G} D(^xf)(x_1)\, dx + R,
$$
where $R$ is supported in $[1,1,1]$.
Further, $D$ is a differential which we can assume to be $G$-equivariant.
Then $D(f)(x_1)$ is of the form
$$
D_1(m) f(x_0. am)|_{a=0},
$$
where $D_1(m)$ is a differential operator in the variable $a$.
Since $D$ is $G$-equivariant, we may replace $f$ with $^{a_0}f$ for some $a_0\in A$.
Since $x_1a_0=x_1$ we get that the above is the same as
$$
D_1(m) f(x_0. aa_0m)|_{a=0}.
$$
This implies that $D_1$ must be of order zero and so $\CT$ is of order zero. 
Restricted to the orbit $x_1.G\cong AM\bs G$ the distribution $\CT$ is given by an integral of the form $\int_{AM\bs G}\phi(y)f(x_1y)\, dy$. (Note that we use the notation without the dot again.)
Invariance implies that $\phi$ is constant.
If $\CT$ is non-zero, then $y\mapsto f(x_1y)$ must be left invariant under $AM$, which implies $\la_1+\la_2-\la_3+\rho=0$, a case we have excluded.
So $\CT$ must be zero.
This shows that any invariant distribution which is zero on the open orbit also vanishes on $x_1.G$.
The remaining orbits are dealt with in a similar fashion.
To prove the Theorem, it remains to show the existence of invariant distributions in the case of unitary parameters.
For this we change our point of view and consider sections of $L_\la$ no longer as functions on $Y$, but as functions on $G^3$ whith values in $V_\sigma=V_{\sigma_1}\otimes V_{\sigma_2}\otimes V_{\sigma_3}$ which spit out $a^{\la+\rho}\sigma(m)$ on the left.
We induce this in the notation by writing $f(x_0y)$ instead of $f(x_0.y)$.
On a given orbit of maximal dimension, there is a standard invariant distribution which, by the lemma, is unique up to scalars and given by
$$
\CT_{\rm st}f\=\al\(\int_Gf(x_0y)\, dy\),
$$
where $\al$ is a linear functional on the space of $M_o$-invariants.
In order to show that this extends to a distribution on $Y$, we 
need to show that the defining integral converges for all $f\in\Ga^\infty(E_{\sigma,\la})$.
This integral equals
\begin{eqnarray*}
\int_G f(x_0y)\, dy &=& \int_Gf(w_0n_0y,w_0y,y)\, dy\\
&=& \int_G\ua(w_0n_0y)^{\rho+\la_1}\ua(w_0y)^{\rho+\la_2}\ua(y)^{\rho+\la_3} f(\uk(x_0y))\, dy.
\end{eqnarray*}
Since $f$ is bounded on $K^3$, it suffices to show the following lemma.

\begin{lemma}\label{2.3}
Let $G=\SO(d,1)^0$ and let $k_0=\uk(w_0n_0)$. Then
$$
\int_G \ua(w_0n_0y)^\rho\ua(w_0y)^\rho\ua(y)^\rho\, dy\ <\ \infty.
$$
\end{lemma}

\begin{conjecture}
The assertion of the lemma should hold for any semisimple group $G$ with finite center and $n_0\in N$ generic.
\end{conjecture}

The conjecture would imply that if $G$ is not locally a product of hyperbolic groups, then the space of invariant trilinear forms on principal series representations is infinite dimensional.

{\bf Proof of Lemma \ref{2.3}:}
Replace the integral over $G$ by an integral over $ANK$ using the Iwasawa decomposition.
Since $\ua(xk)=\ua(x)$ for $x\in G$ and $k\in K$, the $K$-factor is irrelevant and we have to show that
$$
\int_{AN} \ua(k_0an)^\rho\ua(w_0an)^\rho a^\rho\, da\, dn\ <\ \infty.
$$
Now write $w_0n_0=a'n'k_0$. Then $k_0an=(a'n')^{-1}w_0n_0an$ and therefore $\ua(k_0an)=(a')^{-1}\ua(w_0n_0an)$, so it suffices to show
$$
\int_{AN} \ua(w_0n_0an)^\rho\ua(w_0an)^\rho a^\rho\, da\, dn\ <\ \infty.
$$
Next note that $w_0a=a^{-1}w_0$ and so we have $\ua(w_0an)^\rho=a^{-\rho}\ua(w_0n)^\rho$ as well as $\ua(w_0n_0an)^\rho=a^{-\rho}\ua(w_0n_0^an)^\rho$, where $n_0^a=a^{-1}n_0a$.
We need to show
$$
\int_{AN}\ua(w_0n_0^an)^\rho\ua(w_0n)^\rho a^{-\rho}\, da\, dn\ <\ \infty.
$$

This is the point where we have to make things more concrete.
Let $J$ be the diagonal $(d+1)\times (d+1)$-matrix with diagonal entries $(1,\dots,1,-1)$.
Then $\SO(d,1)$ is the group of real matrices $g$ with $g^tJg=J$.
Writing $g=\left(\begin{array}{cc}A & b \\c & d\end{array}\right)$ with $A\in\Mat_d(\R)$, this amounts to
\begin{eqnarray*}
A^tA-c^tc &=& 1\\
A^tb-c^td &=& 0\\
d^2 -b^tb &=& 1
\end{eqnarray*}
The connected component $\SO(d,1)^0$ consists of all matrices $g$ as above with $d>0$.
The maximal compact subgroup $K$ can be chosen to be $\left(\begin{array}{c|c}\SO(d) &  \\\hline  & 1\end{array}\right)$
and $M$ as $\left(\begin{array}{cc|c}\SO(d-1) &  & 0 \\  & 1 & 0 \\\hline 0 & 0 & 1\end{array}\right)$.
Further, we can choose $A$ and $N$ as follows,
\begin{eqnarray*}
A&=& \left\{ \left(\begin{array}{cc|c}1 &  &  \\  & \al & \beta \\\hline  & \beta & \al\end{array}\right): \al>0,\ \al^2-\beta^2=1\right\},\\
&&\\
N &=& \left\{n(x)=\left(\begin{array}{cc|c}1 & -x & x \\ x^t & 1-\frac{|x|^2}2 & |x|^2/2 \\\hline x^t & -|x|^2/2 & 1+\frac{|x|^2}2\end{array}\right): x\in\R^{d-1}\right\},
\end{eqnarray*}
 where we have written $|x|^2=x_1^2+x_2^2+\dots+x_{d-1}^2$.
 Note that the Lie algebra of $A$ is generated by $H=\left(\begin{array}{ccc}1 &  &  \\ & 0 & 1 \\ & 1 & 0\end{array}\right)$.
One derives an explicit formula for the $ANK$-decomposition.
In particular, if $g=\left(\begin{array}{cc}A & b \\c & d\end{array}\right)$ and $\ua(g) =\left(\begin{array}{cc|c}1 &  &  \\  & \al & \beta \\\hline  & \beta & \al\end{array}\right)$, then
$$
\ua(g)^\rho\= (\al+\beta)^{\frac{d-1}2}\= \left(\frac{(d+b_d)}{1+b_1^2+\dots+b_{d-1}^2}\right)^{\frac{d-1}2}.
$$
The Weyl element representative can be chosen to be
$$
w_0\= \left(\begin{array}{ccccc|c}1 &  &  &  &  &  \\ & \ddots &  &  &  &  \\ &  & 1 &  &  &  \\ &  &  & -1 &  &  \\ &  &  &  & -1 &  \\\hline &  &  &  &  & 1\end{array}\right).
$$
So that with $n=n(x)$ for $x\in\R^{d-1}$ we have
\begin{eqnarray*}
\ua(w_0n)^\rho &=& \ua\left(\begin{array}{cc|c}* & * & x_1 \\ * & * & \vdots \\ * & * & x_{d-2} \\ * & * & -x_{d-1} \\ * & * & -\frac{|x|^2}2 \\\hline * & * & 1+\frac{|x|^2}2\end{array}\right)^\rho\\
&=& \left(\frac 1{1+x_1^2+\dots+x_{d-1}^2}\right)^{\frac{d-1}2}
\end{eqnarray*}
We choose $n_0=n(1,0,\dots,0)$ and get with $a=\exp(tH)$,
$$
\ua(w_0n_0^an)^\rho\= \left(\frac 1{1+(e^{-t}+x_1)^2+x_2^2+\dots+x_{d-1}^2}\right)^{\frac{d-1}2}.
$$
With $n=d-1$ and $a^{-\rho}=e^{-\frac{n}2t}$ our assertion boils down to
$$
\int_{\R^n}\int_\R e^{-\frac n2t}(1+|x|^2)^{-\frac n2}(1+|e^{-t}v_1 +x|^2)^{-\frac n2}dt\,dx\ <\ \infty.
$$
Consider first the case $n=1$ and the integral over $x<0$:
$$\int_\R\int_{-\infty}^0 e^{-\frac 12 t}(1+x^2)^{-\frac 12}(1+(e^{-t}+x)^2)^{-\frac 12}dx\, dt
$$
\begin{eqnarray*}
&=& \int_\R\int_0^\infty e^{-\frac 12t}(1+x^2)^{-\frac 12}(1+(x-e^{-t})^2)^{-\frac 12}dx\,dt\\
&=&
\int_\R\int_{-e^{-t}}^\infty e^{-\frac 12t}(1+(x+e^{-t})^2)^{-\frac 12}(1+(x)^2)^{-\frac 12}dx\,dt.
\end{eqnarray*}
Thus it suffices to show the convergence of
$$
\int_{\R}\int_0^\infty e^{-\frac 12t}(1+x^2)^{-\frac 12}(1+(x+e^{-t})^2)^{-\frac 12}dt\,dx.
$$
Setting $y=e^{-t}$ we see that this integral equals
$$
\int_0^\infty\int_0^\infty y^{-1/2}(1+x^2)^{-1/2}(1+(x+y)^2)^{-1/2}dydx
$$
$$
= \int_0^\infty \int_x^\infty (y-x)^{-1/2} (1+x^2)^{-1/2}(1+y^2)^{-1/2} dydx
$$
Since the mapping
$$
x\mapsto \int_x^\infty (y-x)^{-1/2} (1+x^2)^{-1/2}(1+y^2)^{-1/2} dy
$$
is continuous, the integral over $0<x<1$ converges.
It remains to show the convergence of
$$
\int_1^\infty \int_x^\infty (y-x)^{-1/2} (1+x^2)^{-1/2}(1+y^2)^{-1/2} dydx
$$
\begin{eqnarray*}
&=& \int_1^\infty \int_1^\infty x^{1/2}(v-1)^{-1/2} (1+x^2)^{-1/2}(1+v^2x^2)^{-1/2} dvdx\\
&\le& 3^{-1/2} \int_1^\infty \int_1^\infty x^{1/2}(v-1)^{-1/2} (1+x^2)^{-1}(1+v^2)^{-1/2} dvdx\ <\ \infty.
\end{eqnarray*}
here we have used the substitution $y=vx$ and the fact that for $a,b\ge 1$ one has $(1+a)(1+b)\le 3(1+ab)$.

Now for the case $n>1$.
Using polar co-ordinates we compute
$$
\int_{\R^{n-1}}\int_\R\int_\R e^{-tn/2}(1+x^2+|x_r|^2)^{-n/2}(1+)e^{-t}+x)^2+|x_r|^2)^{-n/2} dt\,dx\,dx_r
$$
$$
\ =\ C \int_0^\infty \int_\R\int_\R 
e^{-tn/2}r^{n-2}(1+x^2+r^2)^{-n/2}(1+)e^{-t}+x)^2+r^2)^{-n/2} dt\,dx\,dx_r.
$$
As above we can restrict to the case $x>1$. We get
\begin{eqnarray*}
&=& \int_0^\infty\int_0^\infty\int_0^\infty y^{\frac n2 -1}r^{n-2} (1+x^2+r^2)^{-n/2}(1+(x+y)^2+r^2)^{-n/2}dydxdr\\
&=& \int_0^\infty\int_0^\infty\int_x^\infty (y-x)^{\frac n2 -1}r^{n-2} (1+x^2+r^2)^{-n/2}(1+y^2+r^2)^{-n/2}dydxdr,
\end{eqnarray*}
which equals
$$
\int_0^\infty\int_0^\infty\int_1^\infty (v-1)^{\frac n2 -1}x^{n/2}r^{n-2} (1+x^2+r^2)^{-n/2}(1+v^2x^2+r^2)^{-n/2}dydxdr.
$$
As above, it suffices to restrict the integration to the domain $x>1$.
So one considers
$$
\int_0^\infty\int_1^\infty\int_1^\infty (v-1)^{\frac n2 -1}x^{n/2}r^{n-2} (1+x^2+r^2)^{-n/2}(1+v^2x^2+r^2)^{-n/2}dydxdr.
$$
Choose $0<\eps<1/2$ and write 
\begin{eqnarray*}
(1+x^2+r^2)^{-n/2}&=&(1+x^2+r^2)^{\eps-n/2}(1+x^2+r^2)^{-\eps}\\
&\le& (1+r^2)^{\eps-n/2}(1+x^2)^{-\eps}.
\end{eqnarray*}
So our integral is less than or equal to
$$
\int_0^\infty r^{n-2}(1+r^2)^{\eps-n/2}dr 
$$ 
times
$$
\int_1^\infty\int_1^\infty (v-1)^{\frac n2 -1} x^{n/2}(1+x^2)^{-\eps}(1+v^2x^2)^{-n/2}dvdx
$$
$$
\le C\int_1^\infty\int_1^\infty (v-1)^{\frac n2 -1} x^{n/2}(1+x^2)^{-\eps}(1+v^2)^{-n/2}(1+x^2)^{-n/2}dvdx\ <\ \infty.
$$
\qed

\section{An explicit formula}
Let $d\ge 2$ and $G=\SO(d,1)^0$ if $d>2$.
For $d=2$ let $G$ be the double cover $\tilde\SO(2,1)^o\cong\PGL_2(\R)$.
Then $K=\SO(d)$, $M=\SO(d-1)$ for $d>2$. For $d=2$ we have $K\cong O(2)$ and $M\cong\Z /2\Z$ and in all cases we have $M\bs K\cong S^{d-1}$, the $d-1$ dimensional sphere.
For each $\la\in\a^*$ let $e_\la$ be the class one vector in the associated principal series representation $\pi_\la$ given by $e_\la(ank)=a^{\la+\rho}$.
Let $\la, \mu,\nu\in\a^*$ be imaginary.
Let $\CT_{\rm st}$ be the invariant distribution on $L_{(\la,\mu,\nu)}$ considered
 in the last section.
 We are interested in the growth of $\CT_{\rm st}(e_\la,e_{\mu},e_\nu)$ as a function in $\la$.
First note that the Killing form induces a norm $|.|$ on $\a^*$.

We write $\CT_{\rm st}(\la,\mu,\nu)$ for $\CT_{\rm st}(e_\la,e_{\mu},e_\nu)$ and identifying $\a_\R$ to $\R$ via $\la \mapsto \la(H_0)$ we consider $\CT_{\rm st}$ as a function on $(i\R)^3$.

In this section we will prove the following theorem.

\begin{theorem}\label{3.1} For $\la,\mu,\nu$ imaginary,
$\CT_{\rm st}(\la,{\mu},\nu)$ equals a positive constant times
$$
\frac{\Ga\(\frac{2(\la+\mu-\nu)+n}4\)
\Ga\(\frac{2(\la-\mu+\nu)+n}4\)
\Ga\(\frac{2(-\la+\mu+\nu)+n}4\)
\Ga\(\frac{2(\la+\mu+\nu)+n}4\)}{\Ga\(\frac{2\la+n}2\)
\Ga\(\frac{2\mu+n}2\)
\Ga\(\frac{2\nu+n}2\)},
$$
where $n=d-1$.
So in particular, for fixed imaginary $\mu$ and $\nu$. Then, as $|\la|$ tends to infinity, while $\la$ is imaginary, we have the asymptotic
$$
|\CT_{\rm st}(\la,{\mu},\nu)|\= c\exp\(-\frac \pi 2 |\la|\) |\la|^{\frac d2-2}\( 1+O\(\frac 1{|\la|}\)\),
$$ 
for some constant $c>0$.
\end{theorem}

\prf
The asymptotic formula follows from the explicit expression and the well known asymptotical formula,
$$
|\Ga(\sigma+it)|\= \sqrt{2\pi} \exp\( -\frac\pi 2|t|\)|t|^{\sigma-\frac 12}\( 1+O\(\frac 1{|t|}\)\)
$$
as $|t|\to\infty$, where the real part $\sigma$ is fixed.

Now for the proof of the first assertion.
The permutation group $S_3$ in three letters acts on $(i\R)^3$ by permuting the co-ordinates.
We claim that $\CT_{\rm st}$ is invariant under that action.
So let $\sigma\in S_3$ and let $f=f_1\otimes f_2\otimes f_3$, then
\begin{eqnarray*}
\CT_{\rm st}(\sigma(f)) &=& \int_G f(\sigma^{-1}(x_0.y))\, dy\\
&=& \int_G f(\sigma^{-1}(x_0).y)\, dy
\end{eqnarray*}
Since the open orbit is unique, there is $y_\sigma\in G$ with $\sigma^{-1}(x_0)=x_0.y_\sigma$, and so
\begin{eqnarray*}
\CT_{\rm st}(\sigma(f)) &=& \int_G f(x_0.y_\sigma y)\, dy\\
&=& \int_G f(x_0.y)\, dy\= \CT_{\rm st}(f).
\end{eqnarray*}
So in particular $\CT_{\rm st}(\la,\mu,\nu)$ is invariant under permutations of $(\la,\mu,\nu)$.
 
We have
$$
\CT_{\rm st}(e_\la,e_{\mu},e_\nu)\=\int_G\ua(w_0n_0y)^{\la+\rho}\ua(w_0y)^{\mu+\rho}\ua(y)^{\nu+\rho}\, dy.
$$

Replace the integral over $G$ by an integral over $ANK$ using the Iwasawa decomposition.
Since $\ua(xk)=\ua(x)$ for $x\in G$ and $k\in K$, the $K$-factor is irrelevant and we have to compute
$$
\int_{AN}\ua(w_0n_0an)^{\la+\rho}\ua(w_0an)^{\mu+\rho}a^{\nu+\rho}\, dadn.
$$

Note that $w_0a=a^{-1}w_0$ and so we have $\ua(w_0an)^\rho=a^{-\rho}\ua(w_0n)$ as well as $\ua(w_0n_0an)=a^{-\rho}\ua(w_0n_0^an)$, where $n_0^a=a^{-1}n_0a$.
So the integral equals
$$
\int_{AN}a^{\nu-\la-\mu-\rho}\ua(w_0n_0^an)^{\la+\rho}\ua(w_0n)^{\mu+\rho} \, da\, dn.
$$

Using the explicit approach of the last section we see that $\CT_{\rm st}(\la,\mu,\nu)$ equals
$$
\int_\R\int_{\R^{d-1}}e^{t\left(\nu-\la-\mu-\frac{d-1}2\right)}
(1+|e^{-t}v_1+x|^2)^{-\la-\frac{d-1}2} (1+|x|^2)^{-\mu-\frac{d-1}2}\, dx\,dt,
$$
where $v_1=(1,0,\dots,0)$.

Set $n=d-1=1,2,3,\dots$ and use polar co-ordinates to compute for $n>1$, that
$$
\int_\R\int_{\R^n} e^{t(\nu-\la-\mu-\frac n2)}(1+|e^{-t}v_1+x|^2)^{-\la-\frac n2}(1+|x|^2)^{-\mu-\frac n2}dx\, dt
$$
equals $(n-1)\vol(B_{n-1})$ times
$$
\int_\R\int_{\R}\int_0^\infty r^{n-2} e^{t(\nu-\la-\mu-\frac n2)}(1+(e^{-t}+x)^2+r^2)^{-\la-\frac n2}(1+x^2+r^2)^{-\mu-\frac n2}dr\,dx\, dt.
$$
Define $I_n(a,b,c)$ to be
$$
\int_\R\int_{\R}\int_0^\infty r^{n-2} e^{t(c-a-b)}(1+(e^{-t}+x)^2+r^2)^{-a}(1+x^2+r^2)^{-b}dr\,dx\, dt,
$$
and $I_1(a,b,c)$ to be equal to
$$
\int_\R\int_{\R}  e^{t(c-a-b)}(1+(e^{-t}+x)^2)^{-a}(1+x^2)^{-b}dr\,dx\, dt.
$$
Then, for all $n\in\N$,
$$
\CT_{\rm st}(\la,\mu,\nu)\= c_n I_n\left(\la+\frac n2,\mu+\frac n2,\nu+\frac n2\right),
$$
where $c_1=1$ and $c_n=(n-1)\vol(B_{n-1})$ for $n>1$.
So in particular, $I_n$ is invariant under permutations of the arguments.

Note that
$$
\frac{\partial}{\partial r}\left[ (1+(e^{-t}+x)^2+r^2)^{-a}(1+x^2+r^2)^{-b}\right]
$$
equals
$$
-2ar (1+(e^{-t}+x)^2+r^2)^{-a-1}(1+x^2+r^2)^{-b}\qquad
$$ $$
\qquad -2br (1+(e^{-t}+x)^2+r^2)^{-a}(1+x^2+r^2)^{-b-1},
$$
and $r^{n-2}=\frac{\partial}{\partial r}\left[\frac{r^{n-1}}{n-1}\right]$.
So, integrating by parts, for $n>1$ we compute
\begin{equation}\label{1}
I_n(a,b,c)\= \frac{2a}{n-1}I_{n+2}(a+1,b,c+1)+\frac{2b}{n-1}I_{n+2}(a,b+1,c+1).
\end{equation}
To get a similar result for $I_1$ note that
$$
(1+(e^{-t}+x)^2)^{-a}(1+x^2)^{-b}
$$
equals
$$
-(1+(e^{-t}+x)^2+r^2)^{-a}(1+x^2+r^2)^{-b}|_{r=0}^{r=\infty}
$$ $$
= -\int_0^\infty \frac{\partial}{\partial r}\left[ (1+(e^{-t}+x)^2+r^2)^{-a}(1+x^2+r^2)^{-b}\right]\, dr.
$$
This implies
\begin{equation}\label{2}
I_1(a,b,c)\= 2aI_3(a+1,b,c+1) + 2b I_3(a,b+1,c+1).
\end{equation}
Let $I_n'(a,b,c)$ be the same as $I_n(a,b,c)$ except that there is a factor $x$ in the integrand, i.e.,
$$
\int_\R\int_{\R}\int_0^\infty r^{n-2} e^{t(c-a-b)}(1+(e^{-t}+x)^2+r^2)^{-a}(1+x^2+r^2)^{-b}\,x\,dr\,dx\, dt.
$$
Replacing $x$ with $-x$ first and then with $x+e^{-t}$ yields
\begin{equation}\label{3}
I_n'(a,b,c)\= -I_n'(b,a,c)-I_n(b,a,c-1).
\end{equation}
The last equation also holds for $n=1$.
Integration by parts gives
$$
\int_\R e^{t(c-a-b)}(1+(e^{-t}+x)^2+r^2)^{-a}\, dt\qquad
$$ $$
\qquad = \frac{2a}{a+b-c}\int_\R e^{t(c-a-b)}(1+(e^{-t}+x)^2+r^2)^{-a-1}(e^{-2t}+xe^{-t})\, dt.
$$
So,
$$
I_n(a,b,c)\= \frac{2a}{a+b-c} (I_n(a+1,b,c-1)+I_n'(a+1,b,c)).
$$
Using (\ref{3}) this yields
$$
I_n(a,b,c)=\frac{2a}{c-a-b}I_n'(b,a+1,c),
$$
or
\begin{equation}
I_n'(x,y,z)\= \frac{z-x-y+1}{2y-2} I_n(y-1,x,z).
\end{equation}
Plugging this into (\ref{3}) one gets
\begin{eqnarray*}
\frac{c-a-b+1}{2b-2}I_n(b-1,a,c) &=& \frac{a+b-c-1}{2a-2} I_n(a-1,b,c) - I_n(b,a,c-1).
\end{eqnarray*}
Replacing $(a,b,c)$ with $(a+1,b+1,c+1)$ this gives
$$
\frac{c-a-b}{2b}I_n(b,a+1,c+1)\qquad\qquad
$$ $$
\qquad\qquad \=\frac{a+b-c}{2a}I_n(a,b+1,c+1)-I_n(b+1,a+1,c),
$$
or, using the invariance of $I_n$,
$$
\frac{c-b-a}{2ab}(aI_n(a+1,b,c+1)+bI_n(a,b+1,c+1))\= -I_n(a+1,b+1,c),
$$
and hence
$$
aI_{n+2}(a+1,b,c+1)+bI_{n+2}(a,b+1,c+1)\= \frac{2ab}{a+b-c}I_{n+2}(a+1,b+1,c).
$$
Using (\ref{1}) and (\ref{2}) one gets for $n>1$,
$$
I_n(a,b,c)\=\frac{4ab}{a+b-c}\frac 1{n-1} I_{n+2}(a+1,b+1,c),
$$
and
$$
I_1(a,b,c)\=\frac{4ab}{a+b-c} I_{3}(a+1,b+1,c).
$$
Finally we arrive at
\begin{equation}\label{5}
I_{n+2}(a+1,b+1,c+1)\= \frac{a+b-c-1}{4ab}(n-1)I_n(a,b,c+1)
\end{equation}
for $n>1$ and
\begin{equation}\label{6}
I_{3}(a+1,b+1,c+1)\= \frac{a+b-c-1}{4ab}I_1(a,b,c+1).
\end{equation}
These formulae amount to
$$
\CT_{\rm st}^{n+2}(\la,\mu,\nu)\= d_n \frac{\la+\mu-\nu+\frac n2 -1}{(2\la+n)(2\mu+n)} \CT_{\rm st}^n(\la,\mu ,\nu+1),
$$
where $d_n=(n-1)c_{n+2}/c_n$ if $n>1$ and $d_1=c_3$.
A calculation using the functional equation of the Gamma-function, $\Ga(z+1)=z\Ga(z)$, shows that the right hand side of the claim in the theorem satisfies the same equation as $n$ is replaced by $n+2$. So
 the claim of the theorem for $n$ implies the same claim for $n+2$.
Note that the formula A.5 in \cite{BR} implies the theorem for $n=1$, where one has to take into account that in \cite{BR} a different normalisation is used.
To get our formula from theirs, one has to replace $\la_j$ by $-2\la_j$ in \cite{BR}. 
To finish the proof of the theorem it therefore remains to show the claim for $n=2$.

To achieve this, we proceed in a fashion similar to \cite{BR}.
First note that the group $G=\SL_2(\C)$ is a double cover of $\SO(3,1)^0$, so we might as well use this group for the computation.
For $\la\in\C$ let $V_\la$ denote the space of all smooth functions $f$ on $\C^2$ with $f(az,aw)=|a|^{-2(\la+1)}f(z,w)$ for every $a\in\C^\times$.
Then $G=\SL_2(\C)$ acts on the space $V_\la$ via $\pi_\la(g)f(z,w)=f((z,w)g)$.
This is the principal series representation with parameter $\la$.
For $\la=1$ there is a $G$-invariant continuous linear functional $\L\colon V_1\to\C$, which is unique up to scalars and is given by
$$
\L(f)\=\int_{S^3}f(x)\, dx,
$$ 
where the integral runs over the standard sphere $S^3\subset\C^2$ with the volume element induced by the standard metric on $\C^2\cong\R^4$.
Since $S^3$ equals $(0,1)K$, where $K={\rm SU}(2)$ is the maximal compact subgroup of $G$, the theory of induced representations shows that this functional is indeed invariant under $G$ and is unique with that property.
If $\la$ is imaginary, then the representation $\pi_\la$ is pre-unitary, the inner product being given by $\sp fg=\L(f\bar g)$.
To compute $\CT_{\rm st}^2(\la_1,\la_2,\la_3)$, we will describe this functional on the space $V_{\la_1}\otimes V_{\la_2}\otimes V_{\la_3}$.
Let $\omega\colon\C^2\to\C$ be given by $\omega(a,b)=a_1b_2-a_2b_1$ the bilinear $G$-invariant on $\C^2$.
Let $K_{\la_1,\la_2,\la_3}(v_1,v_2,v_3)$ be the invariant kernel
$$
|\omega(v_1,v_2)|^{(-\la_1+\la_2+\la_3-1)}
|\omega(v_1,v_3)|^{(\la_1-\la_2+\la_3-1)}
|\omega(v_2,v_3)|^{(\la_1+\la_2-\la_3-1)}.
$$
This kernel is invariant under the diagonal action of $G$ and is homogeneous of degree $2(\la_j-1)$ with respect to the variable $v_j$.
From this point the computations run like in \cite{BR} A.4 and A.5.
The theorem follows.
\qed

\appendix
\section{Automorphic coefficients}
Let $\Ga\subset G$ be a uniform lattice and consider the right regular representation of $G$ on $L^2(\Ga\bs G)$, which decomposes as a direct sum,
$$
L^2(\Ga\bs G)\= \bigoplus_{\pi\in\hat G} L^2(\Ga\bs G)(\pi),
$$
where  the isotypic component $L^2(\Ga\bs G)(\pi)$ is zero for $\pi$ outside a countable set of $\pi\in\hat G$ and is always isomorphic to a sum of finitely many copies of $\pi$.
By the Sobolev lemma one has $L^2(\Ga\bs G)^\infty\subset C^\infty(\Ga\bs G)$.
A function $\ph$ is called \emph{pure}, if $\ph\in L^2(\Ga\bs G)(\pi)$ for some $\pi\in\hat G$.
For three pure smooth $\ph_1,\ph_2,\ph_3$ the integral $\int_{\Ga\bs G}\ph_1(x)\ph_2(x)\ph_3(x)\, dx$ exists.
Moreover, fix $\pi_1,\pi_2,\pi_3\in\hat G$ and fix $G$-equivariant embeddings
$\sigma_j\colon V_{\pi_j}\to L^2(\Ga\bs G)$, then
$$
\CT_{\aut}(v_1,v_2,v_3)\= \int_{\Ga\bs G} \sigma_1(v_1)(x)\sigma_2(v_2)(x)\sigma_3(v_3)(x)\, dx
$$
defines a $G$-invariant continuous trilinear form on $V_{\pi_1}^\infty\times V_{\pi_2}^\infty\times V_{\pi_3}^\infty$.
Let $L^2(\Ga\bs G)^K$ be the space of all $K$-invariant vectors in $L^2(\Ga\bs G)$.
Then there is an orthonormal basis $(\ph_i)_{i\in\N}$ of pure vectors in $L^2(\Ga\bs G)^K$.
For each $i$ let $L_i$ be the $G$-stable closed subspace of $L^2(\Ga\bs G)$ generated by $\ph_i$.
Then the spaces $L_i$ are mutually orthogonal.

For fixed pure normalized $\ph,\ph'\in L^2(\Ga\bs G)^K$ we are interested in the growth of the sequence
$$
c_i\df \CT_\aut(\ph_i,\ph,\ph').
$$
Let $(\la_i,\mu,\nu)$ be the induction parameters of the representations given by $(\ph_i,\ph,\ph')$.
First note that the numbers $c_i$ decay exponentially in $|\la_i|$.
So we define the renormalized sequence
$$
b_i\df |c_i|^2\exp(\pi|\la_i|).
$$
In \cite{BR} it is hown that for $d=2$ one has
$$
\sum_{|\la_i|\le T} b_i\ \le\ C\log T,
$$ 
and in \cite{KS}, Theorem 7.6, that for $d\ge 3$,
$$
\sum_{|\la_i|\le T} b_i\ \le\ C T^{2(d-2)}.
$$
We are going to reprove this result for $d\ge 4$.

\begin{theorem}\label{4.1}
Let $d \ge 4$. There is $C>0$ such that for every $T>0$,
$$
\sum_{|\la_i|\le T} b_i\ \le\ C \, T^{2(d-2)}.
$$
\end{theorem}

\prf
By the uniqueness of trilinear forms there exists a constant $a_i$ such that
$$
\CT_\aut(\ph_i,\ph,\ph')\= a_i \CT_{\rm st}(\ph_i,\ph,\ph').
$$
This constant depends on the embeddings of $\pi_\la,\pi_\nu,\pi_\mu$ into $L^2(\Ga\bs G)$, but we can normalize these embeddings by insisting that the standard class one vectors $e_\la, e_\mu, e_\nu$ are mapped to  $\ph_i,\ph,\ph'$.
Thus the number $a_i$ only depends on $\ph_i,\ph,\ph'$.

Fix two smooth pre-unitary $G$-representations $V,V'$ and let $L,L'$ be their Hilbert completions.
Consider the unitary representation of $G\times G$ on $L\otimes L'$ and let $E=(L\otimes L')^\infty$ be its smooth part.
Denote by $\CH(E)$ the real vector space of all Hermitian forms on $E$ and by $\CH^+(E)$ the cone of non-negative forms.

Let $W$ be a smooth pre-unitary admissible representation of $G$.
A $G$-invariant functional $\CT\colon W\otimes V\otimes V'\to\C$ defines a $G$-equivariant linear map $l_\CT\colon V\otimes V'\to W^*$ which extends to a $G$-map $l_\CT\colon E\to \bar W$, where $\bar W$ is the space $W$ with complex conjugate linear structure.
The Hermitian form $H_W$ on $W$ induces a form $H_{\bar W}$ on $\bar W$.
Define a Hermitian form $H_\CT$ on $E$ by $H_\CT=l_\CT^*H_{\bar W}$, i.e., for $u,v\in E$ one has $H_\CT(u,v)=H_{\bar W}(l_\CT(u),l_\CT(v))$.
Finally, suppose that $W=V_\la,V,V'$ are all principal series representations and let $\CT_{\rm st}$ be the standard trilinear form then we write $H_\la^\st$ for the Hermitian form induced on $E$.
We also assume that $V_\la,V,V'$ are cuspidal, i.e., that fixed embeddings into $L^{2,\infty}$ are given.
Then the cuspidal trilinear form $\CT$ induces a Hermitian form $H_\la^\aut$ on $E$.

Next consider the space $L^{2,\infty}(\Ga\bs G\times \Ga\bs G)$.
Let $H_\Delta$ denote the Hermitian form on $L^{2,\infty}(\Ga\bs G\times \Ga\bs G)$ given by restriction to the diagonal, i.e.,
$$
H_\Delta(u,v)\=\int_{\Ga\bs G} u(x,x)\overline{v(x,x)}\, dx.
$$

Let $U$ be a Hilbert space with Hermitian form $H$.
For a closed subspace $L$ of $U$ let $Pr_L$ denote the orthogonal projection from $U$ to $L$ and let $H_L$ be the Hermitian form on $U$ given by
$$
H_L(u,v)\= H(Pr_L(u),Pr_L(v)).
$$
The map $L\mapsto H_L$ is additive and monotonic, i.e., $H_{L+L'}=H_L+H_{L'}$ if $L$ and $L'$ are orthogonal and $H_L\le H_{L'}$ if $L\subset L'$.

\begin{lemma}\label{4.2}
On the space $E$ one has
$$
\sum_i |a_i|^2 H_{\la_i}^\st\ \le\ H_\Delta.
$$
\end{lemma}

\prf
By the uniqueness it follows $H_{\la_i}^\aut= |a_i|^2 H_{\la_i}^\st$ and since all the spaces $L_i$ are orthogonal one also deduces $\sum_i H_{\la_i}^\aut \le H_\Delta$.
\qed

The group $G\times G$ acts on $L^{2,\infty}(\Ga\bs G\times \Ga\bs G)$ and thus on the real vector space of Hermitian forms on $L^{2,\infty}(\Ga\bs G\times \Ga\bs G)$.
By integration, this action extends to a representation of the convolution algebra $C_c^\infty(G\times G)$.
Note that a non-negative function $h\in C_c^\infty(G\times G)$ will preserve to cone of non-negative Hermitian forms.

\begin{lemma}
Let $h\in C_c^\infty(G\times G)$ be non-negative.
Then there exists $C>0$ such that
$$
h.H_\Delta\ \le\ C H_{\aut},
$$
where $H_{\aut}$ is the $L^2$-Hermitian form on $L^{2,\infty}(\Ga\bs G\times \Ga\bs G)$.
\end{lemma}

\prf
The form $h.H_\Delta$ is given by integration over a smooth measure on $\Ga\bs G\times \Ga\bs G$.
Being smooth, this measure is bounded by a multiple of the invariant measure.
\qed

A \emph{positive functional} on $\CH(E)$ is an additive map $\rho\colon \CH^+(E)\to \R_{\ge 0}\cup \{\infty\}$.
For example, every vector $u\in E$ gives rise to a positive functional $\delta_u$ defined as
$$
\delta_u(H)\= H(u,u).
$$
Every positive functional is monotonic and homogeneous, i.e.,
$$
\rho(H)\ \le\ \rho(H')\ \ \ {\rm if}\ \ H\le H',
$$
and
$$
\rho(tH)\= t\rho(H)\ \ \ {\rm for}\ \  t>0.
$$
Let $\rho$ be a positive functional and let $h_\rho(\la)=\rho(H_\la^\st)$.
Then Lemma \ref{4.2} implies
$$
\sum_i h_\rho(\la_i)|a_i|^2\ \le\ \rho(H_\Delta).
$$

\begin{proposition}\label{4.4}
There are constants $T_0, C>0$ such that for every $T\ge T_0$ there is a positive functional $\rho_T$ on $\CH(E)$ with
$$
\rho_T(H_\Delta)\ \le\ C T^{\dim(M\bs K)+1},
$$
and with $h_T=h_{\rho_T}$ one has $h_T(\la)\ge \frac 1C$ for every ${|\la|}\le 2T$.
\end{proposition}

\prf
Let $D$ be a relatively compact open subset of $G$ with $g\in D_1\Rightarrow g^{-1}\in D$.
Assume further that $D$ is invariant under right and left translations by elements of $K$.
Let $h\in C_c^\infty(G\times G)$ such that $h$ is constantly one on $D\times D$ and such that $h$ is invariant under left and right translates by elements of $K$.
the space $E$ can be identified with $C^\infty(M\bs K\times M\bs K)$.
Fix a norm $\norm\cdot$ on the Lie algebra of $AN$.

Note that the map
\begin{eqnarray*}
\psi\colon A\times N\times M &\to& P\bs G\times P\bs G\\
(a,n,m)&\mapsto & (w_0anm,w_0n_0anm)
\end{eqnarray*}
is a diffeomorphism onto an open subset of $(P\bs G)^2\cong (M\bs K)^2$.
Let $F=\psi(1\times 1\times M)=$ the $M$-orbit of the point $(w_0,w_0n_0)$.
Let $C_1,T>0$ be large and set
$$
S\= \left\{\psi(a,n,m): \norm{\log an}< \frac 1{C_1T}, m\in M\right\}.
$$
Let $u_T$ be a nonnegative smooth function on $(M\bs K)^2$ with support in $S$ such that
$$
\int u\, dxdy\= 1,\ \ \ \ \int|u|^2\, dxdy\ \le\ T^{\dim(M\bs K)+1}.
$$
It is possible to construct such a function since the codimension of $F$ equals $\dim(M\bs K)+1$.
With these data define a positive functional $\rho_T$ on $\CH(E)$ by
$$
\rho_T(H)\= \delta_{u_T}(h.H).
$$
Then one has
$$
\rho_T(H_\Delta)\=\rho_u(h.H_\Delta)\ \le\ C'\rho_u(H_{st})\ \le\ C T^{\dim(M\bs K)+1}.
$$\qed

Recall the standard trilinear form $\CT_\st$ on $V_\la\otimes V\otimes V'$ and the corresponding linear map $l_{\CT_\st}\colon V\otimes V'\to \bar V_\la$.
We identify the space $\bar V_\la$ with $C^\infty(M\bs K)$ and let $z$ denote the point $M 1\in M\bs K$.
The Dirac measure $\delta_z$ at $z$ is a continuous linear functional on $C^\infty(M\bs K)$.
We get an induced Hermitian form $H_z$ on $E$ defined by
$$
H_z(u,v)\= l_{\CT_\st}(u)(z)\overline{l_{\CT_\st}(v)(z)}.
$$

\begin{lemma}On $E$ one has
$$
H_\la^\st\=\int_K (k,k).H_z\, dk.
$$
\end{lemma}

\prf
This follows from the fact that the invariant Hermitian form on $\bar V_\la$ is equal to $\int_K\pi(k)\tilde H_z\, dk$,
where $\tilde H_z$ is the Hermitian form on $\bar V_\la$ given by $\tilde H_z(u,v)=u(z)\overline {v(z)}$.
\qed

Since the test function $h$ was assumed to be $K\times K$-invariant it follows that
$$
h_T(\la)\=\rho_T(H_\la^\st)\=\delta_u(h.H_z).
$$
So in order to establish a lower bound for $\rho_T(H_\la^\st)$ it suffices to give a lower bound for $\delta_u(x.H_z)=|\sp{\pi(x)f_z}u|^2$ for $x\in D$, where $f_z$ is the linear functional on $E$ given by $f_z(u)=\delta_z(l_{\CT_\st}(u))$.

\begin{lemma}\label{4.6}
There is $C>0$ such that 
for $T_0$ large enough there exists an open, non-empty subset $D_0\subset D$ such that for $T\ge T_0$, $|\la|\le T$, and $x\in D_0$ we have $|\sp{\pi(x)f_z}{u_T}|\ge \frac 1C$.
\end{lemma}

\prf
Note first that with $x=ank\in G$ we have 
$$
\sp{\pi(x)f_z}{u_T}\=\int_{(M\bs K)^2}u_T(k_2,k_3)\phi(k,k_2^{ank},k_3^{ank})\, dk_2dk_3.
$$
We next derive an explicit formula for $\phi$ on the open orbit.
Recall that
\begin{eqnarray*}
\CT_{\rm st}(f) &=& \int_G f(y,w_)y,w_0n_0y)\, dy\\
&=& \int_{ANK} f(ank,w_0ank,w_0n_0ank)\, da\,dn\,dk,
\end{eqnarray*}
and that the last integrand equals
$$
a^{\la+\rho}\ua(w_0an)^{\mu+\rho}\ua(w_0n_0an)^{\nu+\rho} f(k,\uk(w_0an)k,\uk(w_0n_0an)k)\, da\,dn\,dk.
$$
We deduce that
$$
\phi(k,\uk(w_0an)k,\uk(w_0n_0an)k) \= a^{\la+\rho}\ua(w_0an)^{\mu+\rho}\ua(w_0n_0an)^{\nu+\rho} W(an),
$$
where $W$ is the Radon-Nykodim derivative of $dadndk$ against the measure on $(M\bs K)^3$. 
In particular, this function is positive and smooth and does not depend on $\la$.
It follows that the functional $\pi(x)f_z$ is given by integration against the function
$$
\pi(ank)f_z(w_0a'n'm',w_0n_0a'n'm')\= \phi(k,w_0aa'n''m'k,w_0n_0aa'n''m'k),
$$
where $n''=a^{-1}n'a\, m'n(m')^{-1}$.
Note that $Mk=Mm'k$ and so the above equals
$$
(aa')^{\la+\rho}\ua(w_0aa'n'')^{\mu+\rho}\ua(w_0n_0aa'n'')^{\nu+\rho} W(aa'n'').
$$
Note in particular, that the function $\pi(x)f_z$ is invariant under $M$.

The set $S$ is a tubular neighbourhood of the compact set $F$.
We have natural local coordinates with $F$-directions, which can be identified with open parts of $M$-orbits and $F^\perp$-directions, which can be viewed as $AN$-orbits.
The above formula shows that the gradient of $\pi(x)f_z$ is zero along the $M$-orbits.
Along the $AN$-orbits this gradient is bounded by a constant times $T$, if we leave $\mu$ and $\nu$ fixed and restrict to $|\la|\le T$.
Let $C_0>0$ be large and let $D_0\subset D$ be the subset of elements $x\in D$ such that $|\pi(x)f_z|< C_0$.
This set is open and nonempty if $C_0$ is large enpough and it does not depend on $\la$.

Lemma \ref{4.6} implies Proposition \ref{4.4} and the latter implies Theorem \ref{4.1} as follows.

Consider the inequality $\sum_i|a_i|^2\rho_T(H_{\la_i}^\st)\le\rho_T(H_\Delta)$.
By Proposition \ref{4.4} the right hand side is bounded by $C T^{2\dim (M\bs K)}$ and since $\dim M\bs K =\dim P\bs G =\dim N= d-1$ we get
$$
\sum_i|a_i|^2\rho_T(H_{\la_i}^\st)\ \le\ CT^{d}.
$$
Restricting the sum to those $i$ with $|\la_i|\le 2T$ and using the second assertion of Proposition \ref{4.4} we get
$$
\sum_{|\la_i|\le T}|a_i|^2\ \le \ C T^{d}.
$$
According to Theorem \ref{3.1} we have $b_i|\la_i|^{4-d}\sim |a_i|^2$ and so
$$
\sum_{|\la_i|\le T}b_i|\la_i|^{4-d}\ \le \ C T^{d}.
$$
This implies the theorem.
\qed

{\small Mathematisches Institut\\
Auf der Morgenstelle 10\\
72076 T\"ubingen\\
Germany\\
\tt deitmar@uni-tuebingen.de}

\end{document}